\documentstyle[12pt]{article}
\newlength{\abstractwidth}
\renewcommand{\thefootnote}{\fnsymbol{footnote}}
\renewcommand{\thanks}[1]{\footnote{#1}} % Use this for footnotes
\newcommand{\starttext}{ \setcounter{footnote}{0}
\renewcommand{\thefootnote}{\arabic{footnote}}}

\newcommand{\be}{\begin{equation}}
\newcommand{\bea}{\begin{eqnarray}}
\newcommand{\eea}{\end{eqnarray}} \newcommand{\ee}{\end{equation}}
 \newcommand{\<}{\langle}
\renewcommand{\>}{\rangle} \def\ba{\begin{eqnarray}}
\def\ea{\end{eqnarray}}

\def\o{\omega}
\def\Re{{\rm Re}}

\def\log{\,{\rm log}\,}

\def\o{\omega}

\def\o{\omega}

\def\na{\nabla}

\def\p{\partial}

\def\ddb{{\partial\bar\partial}}

\def\na{{\nabla}}

\def\[{{\bf [}}
\def\]{{\bf ]}}

\begin{document}
\starttext \baselineskip=18pt \setcounter{footnote}{0}
\newtheorem{theorem}{Theorem}
\newtheorem{lemma}{Lemma}
\newtheorem{corollary}{Corollary}
\newtheorem{definition}{Definition}
\newtheorem{conjecture}{Conjecture}
\newtheorem{proposition}{Proposition}

\begin{center}
{\Large \bf UNIFICATION OF THE K\"AHLER-RICCI AND ANOMALY FLOWS
\footnote{Work supported in part by the National Science Foundation under grant DMS-12-66033.}}

\medskip
\centerline{ Teng Fei and Duong H. Phong}

\medskip

\begin{abstract}

A new formulation of the Anomaly flow in the case of vanishing slope parameter is given, where the dependence on the global section of the canonical bundle appears only in the initial data. This allows a natural unification of the Anomaly flow with the K\"ahler-Ricci flow.

\end{abstract}

\end{center}

\baselineskip=15pt
\setcounter{equation}{0}
\setcounter{footnote}{0}

\section{Introduction}
\setcounter{equation}{0}

The idea of using a geometric flow to implement a cohomological constraint on a metric in the absence of an analogue of the $\partial\bar\partial$ lemma was introduced in \cite{PPZ1, PPZ2}. The specific case of the conformally balanced condition arising from supersymmetric compactifications of the heterotic string was considered there and generalized further in \cite{PPZ5}. Many other conditions and flows have been introduced since, including dual Anomaly flows \cite{FP} and flows motivated by Type II A and Type II B string compactifications in \cite{P, FPPZ}. Anomaly flows appear to be a flexible and powerful method, as they have led to new proofs of major results in geometry such as Yau's theorem \cite{Y} on the existence of K\"ahler Ricci-flat metric and the Fu-Yau solution \cite{FY1,FY2} of the Hull-Strominger system
\cite{PPZ3, PPZ4, FHP}.

\smallskip
A flow is usually given by a vector field on the configuration space and the prescription of an initial data. In the Anomaly flows considered in \cite{PPZ1, PPZ2, PPZ4}, the underlying manifold $X$ is complex and equipped with a non-vanishing top holomorphic form $\Omega$. The form $\Omega$ appears explicitly in the vector field on the space of Hermitian metrics defining the flow (see \cite{PPZ2}, eq.(1.9)). This explicit appearance of $\Omega$ seems to set Anomaly flows apart from more familiar flows such as the K\"ahler-Ricci flow, and prevent the direct use of many powerful techniques which had been developed for these flows.

\smallskip
The main purpose of the present note is to show that, in the simpler case with parameter $\alpha'=0$, the dependence of the vector field of the Anomaly flow on $\Omega$ can be eliminated by a suitable rescaling of the evolving metric:

\begin{theorem}
\label{main}
Let $X$ be a complex manifold of dimension $m\geq 2$ equipped with a nowhere holomorphic $(m,0)$-form $\Omega$. Assume that $t\to\omega(t)$ is a flow of Hermitian metrics satisfying
\bea
\label{anomaly}
\p_t(\|\Omega\|_\o\o^{m-1})=i\p\bar\p\o^{m-2}
\eea
and $d(\|\Omega\|_{\o(t)}\o(t)^{m-1})=0$ for each $t$. Set for each $t$
\bea
\eta(t)=\|\Omega\|_{\o(t)}\o(t).
\eea
Then the Hermitian metrics $\eta(t)$ satisfy the conformally balanced condition $d(\|\Omega\|_\eta^2\eta^{m-1})=0$, and they evolve according to
\bea
\label{eta-flow}
\p_t\eta(t)=-{1\over m-1}(\tilde R_{\bar kj}(\eta)+{1\over 2}(T\circ\bar T)_{\bar kj}).
\eea
Here $\tilde R_{\bar kj}:=R^p{}_{p\bar kj}=-g^{p\bar q}g_{\bar km}\p_{\bar q}(g^{m\bar\ell}\p_pg_{\bar\ell j})$ is the Chern-Ricci tensor of $\eta:=ig_{\bar kj}dz^j\wedge d\bar z^k$ and $(T\circ\bar T)_{\bar kj}:=T_{\bar kpq}\bar T_j{}^{pq}$, where $T=i\p\eta=
{1\over 2}T_{\bar kpq}dz^q\wedge dz^p\wedge d\bar z^k$ is the torsion of $\eta$.
\end{theorem}

The form $\Omega$ has cancelled out from the vector field $\p_t\eta$, as desired. From the point of view of $\eta(t)$, the only dependence on $\Omega$ of the Anomaly flow resides now in the conformally balanced condition for the initial data $\eta(0)$. Thus the flow defined by the right hand side of (\ref{eta-flow}) with an arbitrary initial Hermitian metric can be viewed as a generalization of the Anomaly flow with $\alpha'=0$ to arbitrary complex manifolds $X$. When $X$ is compact, it is not difficult to see, as we shall show in detail later, that the flow
(\ref{eta-flow}) preserves the K\"ahler property and reduces to the K\"ahler-Ricci flow if the initial data is K\"ahler. We can then formulate the following theorem, which is essentially Theorem 1 combined with the uniqueness of solutions of parabolic flows on compact manifolds, and which unifies the K\"ahler-Ricci flow with the Anomaly flow:

\begin{theorem}
\label{unification}
Let $X$ be a compact complex manifold of dimension $m\geq 1$. Consider the flow $t\to\eta(t)$ of Hermitian metrics defined by
\bea
\label{eta-flow1}
\p_t\eta(t)=-(\tilde R_{\bar kj}(\eta)+{1\over 2}(T\circ\bar T)_{\bar kj}).
\eea
with initial data a Hermitian metric $\eta(0)$

{\rm (i)} The flow is parabolic, and for any initial data $\eta(0)$, it admits a unique smooth solution
in some maximal time interval $[0,T)$ with $T>0$.

{\rm (ii)} If $\eta(0)$ is K\"ahler (this includes the general case in dimension $m=1$), then $\eta(t)$ remains K\"ahler for all time $t\in [0,T)$, and the flow reduces to the K\"ahler-Ricci flow,
\bea
\label{kaehler-ricci}
\p_t\eta_{\bar kj}=-R_{\bar kj}(\eta).
\eea

{\rm (iii)} Assume that $m\geq 2$, and $X$ admits a nowhere vanishing holomorphic $(m,0)$-form $\Omega$. If $\eta(0)$ is conformally balanced in the sense that $d(\|\Omega\|_{\eta(0)}^2\eta^{m-1}(0))=0$,  then $\eta(t)$ remains conformally balanced for all time $t\in [0,T)$, and the flow reduces to the Anomaly flow
(\ref{anomaly}), after rescaling $t\to (m-1)t$.

\end{theorem}

Parts (i) and (ii) of Theorem \ref{unification} are elementary, and have been noted by
Streets and Tian \cite{ST} who proposed the family of flows of the form
\bea
\label{eta-ST}
\p_t\eta_{\bar kj}=-(\tilde R_{\bar kj}+Q_{\bar kj}(T,\bar T))
\eea
as generalizations of the K\"ahler-Ricci flow to arbitrary complex manifolds, where $Q(T,\bar T)$ is a $(1,1)$-form which is linear in each factor $T$ and $\bar T$. Among these, Ustinovskiy \cite{U} has identified the same combination $Q={1\over 2}T\circ\bar T$ as in (\ref{eta-flow}) as a flow that preserves the Griffiths positivity and the dual Nakano-positivity of the tangent bundle. In the case of the K\"ahler-Ricci flow, the preservation of the positivity of the bisectional curvature was proven by Bando \cite{Ba} and Mok \cite{Mok} and is a particularly important property of the flow with many applications, see e.g. \cite{PSSW}. Ustinovskiy's result \cite{U} suggests that some generalizations in (\ref{eta-ST}) may be better behaved than others. With Theorem \ref{main}, we see that the Anomaly flow with $\alpha'=0$ also singles out the particular combination $Q={1\over 2}T\circ\bar T$.

\medskip
We note that many other generalizations of the K\"ahler-Ricci flow to the non-K\"ahler setting have been proposed in the literature, including in \cite{G}, \cite{BV}, and \cite{TW}.

\section{Proof of Theorem \ref{main}}
\setcounter{equation}{0}

First, we note that the rescaled metric
$\eta=\|\Omega\|_\o\o$ satisfies
\bea
\|\Omega\|_\o^{2-m}=\|\Omega\|_\eta^2,
\qquad
\|\Omega\|_\o\o^{m-1}=\|\Omega\|_\eta^2\eta^{m-1},
\eea
and hence the Anomaly flow (\ref{anomaly}) can be expressed in terms of $\eta$ as
\bea
\label{eta-flow2}
\p_t(\|\Omega\|_\eta^2\eta^{m-1})
=i\ddb \, (\|\Omega\|_\eta^2\eta^{m-2}).
\eea

\subsection{Elimination of $\p_t\|\Omega\|_\eta^2$}

Carrying out the differentiation in time gives
\bea
\|\Omega\|_\eta^2(\p_t\log\|\Omega\|_\eta^2\eta^{m-1}
+
(m-1)\p_t\eta\wedge\eta^{m-2})=i\p\bar\p(\|\Omega\|_\eta^2\eta^{m-2})
\eea
Let $\Lambda$ be the usual Hodge operator, which is the adjoint of the operator $\phi\to\o\wedge\phi$. Its precise expression in components and normalization can be found in Appendix C. Since
\bea
\p_t\log\|\Omega\|_\eta^2=\p_t\log\eta^{-m}=-g^{j\bar k}\p_tg_{\bar kj}=-(\Lambda \p_t\eta)
\eea
we obtain
\bea
\|\Omega\|_\eta^2(-(\Lambda\p_t\eta)\eta^{m-1}+(m-1)\p_t\eta\wedge\eta^{m-2})=i\p\bar\p(\|\Omega\|_\eta^2\eta^{m-2}).
\eea
We take now the Hodge $\star$ operator of both sides, using the formulas in Appendix B. We find
\bea
\|\Omega\|_\eta^2
\bigg\{-(m-1)!(\Lambda\p_t\eta)\eta
+
(m-1)(m-2)!(-\p_t\eta+(\Lambda\p_t\eta)\eta)\bigg\}=\star i\p\bar\p(\|\Omega\|_\eta^2\eta^{m-2})
\eea
The term $(\Lambda\p_t\eta)$ cancels out from the left hand side, and we obtain the following equation
\bea
(m-1)!\|\Omega\|_\eta^2\p_t\eta=-\star i\p\bar\p(\|\Omega\|_\eta^2\eta^{m-2}).
\eea

\subsection{Elimination of $\|\Omega\|_\eta^2$}

It is now easy to see that the explicit appearance of the term $\|\Omega\|_\eta^2$ can be eliminated from the flow, using the torsion constraints. Indeed, before taking the Hodge $\star$ operator, the right hand side of the flow can be expressed as
\bea
 i\p\bar\p(\|\Omega\|_\eta^2\eta^{m-2})
 &=&i\p(\|\Omega\|_\eta^2\bar\p\log\|\Omega\|_\eta^2 \eta^{m-2}+\|\Omega\|_\eta^2\bar\p\eta^{m-2})
 \nonumber\\
 &=&
 \|\Omega\|_\eta^2
 \bigg\{i\p\log\|\Omega\|_\eta^2 \wedge \bar\p\log\|\Omega\|_\eta^2\wedge\eta^{m-2}+i\p\bar\p\log\|\Omega\|_\eta^2\wedge\eta^{m-2}
 \nonumber\\
 &&
 -
 i\bar\p\log\|\Omega\|_\eta^2\p\eta^{m-2}
 +
 i\p\log\|\Omega\|_\eta^2\bar\p\eta^{m-2}+i\p\bar\p\eta^{m-2}\bigg\}
 \eea
 However, from Lemma 4 in \cite{PPZ2}, with $a=2$, we have
 \bea
 \tau_\ell=g^{j\bar k}T_{\bar kj\ell}=\p_\ell\log\|\Omega\|_\eta^2
 \eea
 and in general, $Ric(\eta)=\ddb\log\|\Omega\|_\eta^2$. Thus the above equation can be rewritten as
 \bea
 i\p\bar\p(\|\Omega\|_\eta^2\eta^{m-2})&=&
  \|\Omega\|_\eta^2
 \bigg\{i\tau\wedge \bar\tau \wedge\eta^{m-2}+iRic(\eta)\wedge \eta^{m-2}+i\p\bar\p\eta^{m-2}
 \nonumber\\
 &&
 -i\bar\tau\wedge\p\eta^{m-2}+i\tau \wedge\bar\p\eta^{m-2}\bigg\}.
 \eea
Returning to the Anomaly flow, it reduces now to the following simpler expression
\bea
(m-1)!\p_t\eta
=
-\star\bigg\{i\tau\wedge \bar\tau \wedge\eta^{m-2}+iRic(\eta)\wedge \eta^{m-2}+i\p\bar\p\eta^{m-2}
 -i\bar\tau\wedge\p\eta^{m-2}+i\tau \wedge\bar\p\eta^{m-2}\bigg\}.
 \nonumber\\
 \eea
 Next, we note that
 \bea
 \p\eta^{m-2}&=&(m-2)\p\eta\wedge\eta^{m-3}=-i(m-2)T\wedge\eta^{m-3}
 \nonumber\\
 \bar\p\eta^{m-2}
 &=& i(m-2)\bar T\wedge \eta^{m-3}
 \eea
 and
\bea
i\ddb \eta^{m-2}&=&(m-2) i\p(\bar\p\eta\wedge \eta^{m-3})
=
(m-2)(i\p\bar\p\eta\wedge \eta^{m-3}-i(m-3)\bar\p\eta\wedge \p\eta\wedge \eta^{m-4})
\nonumber\\
&=&
(m-2)(i\p\bar\p\eta\wedge \eta^{m-3}
-
i(m-3)\bar T\wedge T\wedge \eta^{m-4}).
\eea
Collecting all the terms, we obtain
\bea
(m-1)!\p_t\eta
&=&-\star\bigg\{(i\tau\wedge\bar\tau+iRic(\eta))\wedge \eta^{m-2}
+(m-2)(i\ddb\eta-\bar\tau\wedge T-\tau\wedge \bar T)\wedge\eta^{m-3}
\nonumber\\
&&\quad-(m-2)(m-3)i\bar T\wedge T\wedge \eta^{m-4}\bigg\}
\eea

\subsection{The Hodge $\star$ of the individual terms}

Applying the formulas for the Hodge $\star$ operator given in the appendices, we obtain immediately
\bea
\star[(i\tau\wedge\bar\tau+iRic(\eta))\wedge \eta^{m-2}]
=
(m-2)![-(i\tau\wedge\bar\tau+iRic(\eta))+(\Lambda (i\tau\wedge\bar\tau+iRic(\eta)))\eta]
\nonumber
\eea
and
\bea
&&
\star[(m-2)(i\ddb\eta-\bar\tau\wedge T-\tau\wedge \bar T)\wedge \eta^{m-3}]
=(m-2)![-\Lambda (i\ddb\eta-\bar\tau\wedge T-\tau\wedge \bar T)\nonumber\\
&&
\qquad\qquad\qquad
+{1\over 2}\Lambda^2(i\ddb\eta-\bar\tau\wedge T-\tau\wedge \bar T)\,\eta]
\eea
and
\bea
\star[(m-2)(m-3)iT\wedge\bar T\wedge\eta^{m-4}]
=
(m-2)![-{1\over 2}\Lambda^2(iT\wedge\bar T)+{1\over 6}\Lambda^3(iT\wedge\bar T)\eta]
\eea
The appearance of a common factor $(m-2)!$ in all the terms of the right hand side allows us to cancel this factor, and obtain a generalization of the anomaly flow including to dimension $m=2$, defined by
\bea
-(m-1)\p_t\eta&=&
-(i\tau\wedge\bar\tau+iRic(\eta))+(\Lambda (i\tau\wedge\bar\tau+iRic(\eta)))\eta
\nonumber\\
&&
-\Lambda (i\ddb\eta-\bar\tau\wedge T-\tau\wedge \bar T)+{1\over 2}\Lambda^2(i\ddb\eta-\bar\tau\wedge T-\tau\wedge \bar T)\,\eta
\nonumber\\
&&
-{1\over 2}\Lambda^2(iT\wedge\bar T)+{1\over 6}\Lambda^3(iT\wedge\bar T)\eta
\nonumber\\
&=&A+B\eta
\eea
where we have defined the $(1,1)$-form $A$ and the scalar function $B$ by
\bea
A&=&-(i\tau\wedge\bar\tau+iRic(\eta))-\Lambda (i\ddb\eta-\bar\tau\wedge T-\tau\wedge \bar T)
\nonumber\\
&&-{1\over 2}\Lambda^2(iT\wedge\bar T)
\nonumber\\
B&=&\Lambda (i\tau\wedge\bar\tau+iRic(\eta))+{1\over 2}\Lambda^2(i\ddb\eta-\bar\tau\wedge T-\tau\wedge \bar T)
+{1\over 6}\Lambda^3(iT\wedge\bar T).
\eea

\subsection{Evaluation of $i\ddb\eta$}

We quote from \cite{PPZ2}, eq. (2.52)
\bea
(i\ddb\eta)_{\bar kj\bar \ell m}
=
R_{\bar kj\bar\ell m}
-
R_{\bar km\bar \ell j}+R_{\bar\ell m\bar kj}-R_{\bar\ell j\bar km}-g^{s\bar r}T_{\bar r jm}\bar T_{s\bar k\bar\ell}
\eea
It follows that
\bea
(\Lambda i\ddb\eta)_{\bar\ell m}
=i^{-1}g^{j\bar k}(i\ddb\eta)_{\bar kj\bar\ell m}=
i^{-1}(\tilde R_{\bar\ell m}+R_{\bar\ell m}-g^{s\bar r}g^{j\bar k}T_{\bar r jm}\bar T_{s\bar k\bar\ell})
\eea
or, in terms of forms,
\bea
\Lambda i\ddb\eta
=-i\tilde Ric(\eta)-iRic(\eta)+i(T\bar T)
\eea
where the $(1,1)$-form $T\bar T)$ is defined by
\bea
T\bar T=(T\bar T)_{\bar \ell m}dz^m\wedge d\bar z^\ell,
\qquad
(T\bar T)_{\bar \ell m}=g^{s\bar r}g^{j\bar k}T_{\bar r jm}\bar T_{s\bar k\bar\ell}.
\eea
As a consequence, the Ricci-Chern terms cancel and we obtain
\bea
iRic(\eta)+\Lambda i\ddb\eta=
-i\tilde Ric(\eta)+i(T\bar T).
\eea
Similarly,
\bea
\Lambda^2i\ddb\eta=-2R+|T|^2
\eea
and the scalar curvature cancels between the terms $\Lambda iRic(\eta)$ and $\Lambda^2i\ddb\eta$,
\bea
\Lambda iRic(\eta)+{1\over 2}\Lambda^2i\ddb\eta={1\over 2}|T|^2.
\eea

It is then convenient to isolate torsion and non-torsion terms in the coefficients $A$ and $B$ as follows
\bea
A&=&i\tilde Ric(\eta)-i(T\bar T)-i\tau\wedge\bar\tau+\Lambda(\tau\wedge \bar T+\bar\tau\wedge T)-{1\over 2}\Lambda^2(iT\wedge\bar T)
\nonumber\\
B&=&{1\over 2}|T|^2+|\tau|^2-{1\over 2}\Lambda^2(\bar\tau\wedge T+\tau\wedge\bar T)+{1\over 6}\Lambda^3(iT\wedge\bar T).
\eea

\subsection{Evaluation of $iT\wedge\bar T$, $\Lambda(iT\wedge\bar T)$, $\Lambda^2(iT\wedge\bar T)$}

The components of $iT\wedge\bar T$ can be expressed as, upon antisymmetrization,
\bea
(iT\wedge\bar T)_{\bar kj\bar\beta\alpha\bar\gamma\ell}
&=&-i(T_{\bar kj\ell}\bar T_{\alpha\bar\beta\bar\gamma}-T_{\bar k\alpha\ell}\bar T_{j\bar\beta\bar\gamma}-T_{\bar kj\alpha}\bar T_{\ell\bar\beta\bar\gamma}
-T_{\bar\beta j\ell}\bar T_{\alpha\bar k\bar\gamma}+T_{\bar\beta\alpha\ell}\bar T_{j\bar k\bar\gamma} +T_{\bar\beta j\alpha}\bar T_{\ell\bar k\bar\gamma}
\nonumber\\
&&
-T_{\bar\gamma j\ell}\bar T_{\alpha\bar\beta\bar k}+T_{\bar\gamma\alpha \ell}\bar T_{j\bar\beta\bar k}+T_{\bar\gamma j\alpha}\bar T_{\ell\bar\beta\bar k})
\eea
It follows that
\bea
(\Lambda iT\wedge\bar T)_{\bar kj\bar\beta\alpha}
&=&-g^{\ell\bar\gamma}(T_{\bar kj\ell}\bar T_{\alpha\bar\beta\bar\gamma}-T_{\bar k\alpha\ell}\bar T_{j\bar\beta\bar\gamma})
+g^{\ell\bar\gamma}(-T_{\bar\beta j\ell}\bar T_{\alpha\bar k\bar\gamma}+T_{\bar\beta\alpha\ell}\bar T_{j\bar k\bar\gamma} )+g^{\ell\bar\gamma}
T_{\bar\gamma j\alpha}\bar T_{\ell\bar\beta\bar k}
\nonumber\\
&&
-T_{\bar kj\alpha}\bar\tau_{\bar\beta}+T_{\bar\beta j\alpha}\bar\tau_{\bar k}
-\tau_j\bar T_{\alpha\bar\beta\bar k}+\tau_{\alpha}\bar T_{j\bar\beta\bar k}
\eea
Note that
\bea
\tau\wedge \bar T&=&\tau_\alpha dz^\alpha\wedge {1\over 2}\bar T_{j\bar\beta \bar k}d\bar z^k\wedge d\bar z^\beta\wedge dz^j
\nonumber\\
&=&
{1\over 2^2}(\tau_\alpha \bar T_{j\bar \beta\bar k}-\tau_j\bar T_{\alpha\bar\beta\bar k})
dz^\alpha\wedge d\bar z^\beta\wedge dz^j\wedge d\bar z^k
\eea
and hence
\bea
(\tau\wedge \bar T)_{\bar kj\bar\beta\alpha}=\tau_\alpha \bar T_{j\bar \beta\bar k}-\tau_j\bar T_{\alpha\bar\beta\bar k}.
\eea
It follows that
\bea
(\Lambda \tau\wedge\bar T)_{\bar\beta\alpha}
=i^{-1}g^{j\bar k}(\tau_\alpha \bar T_{j\bar \beta\bar k}-\tau_j\bar T_{\alpha\bar\beta\bar k})=
i\tau_\alpha \bar\tau_{\bar\beta}+ig^{j\bar k}\tau_j\bar T_{\alpha\bar\beta\bar k}.
\eea
and
\bea
\Lambda^2\tau\wedge\bar T
=2|\tau|^2.
\eea
Returning to the earlier identity, we can now compute $\Lambda^2 iT\wedge\bar T$,
\bea
(\Lambda^2iT\wedge\bar T)_{\bar\beta\alpha}
&=&ig^{j\bar k}\bigg\{g^{\ell\bar\gamma}(T_{\bar kj\ell}\bar T_{\alpha\bar\beta\bar\gamma}-T_{\bar k\alpha\ell}\bar T_{j\bar\beta\bar\gamma})
+g^{\ell\bar\gamma}(-T_{\bar\beta j\ell}\bar T_{\alpha\bar k\bar\gamma}+T_{\bar\beta\alpha\ell}\bar T_{j\bar k\bar\gamma} )+g^{\ell\bar\gamma}
T_{\bar\gamma j\alpha}\bar T_{\ell\bar\beta\bar k}\bigg\}
\nonumber\\
&&+(\Lambda(\tau\wedge\bar T+\bar\tau\wedge T))_{\bar\beta\alpha}
\nonumber\\
&=&
ig^{\ell\bar\gamma}(\tau_\ell \bar T_{\alpha\bar\beta\bar\gamma}+\bar\tau_{\bar\gamma}T_{\bar\beta\alpha\ell})
-i (T\circ \bar T)_{\bar\beta\alpha}
-2i(T\bar T)_{\bar\beta\alpha}\nonumber\\
&&+(\Lambda(\tau\wedge\bar T+\bar\tau\wedge T))_{\bar\beta\alpha}
\eea
where we have defined the $(1,1)$-form $T\circ \bar T$ by
\bea
(T\circ\bar T)_{\bar\beta\alpha}
=
g^{\ell\bar\gamma}g^{j\bar k}T_{\bar\beta j\ell}\bar T_{\alpha\bar k\bar\gamma}.
\eea
In intrinsic notation, this can be expressed as
\bea
(\Lambda^2iT\wedge\bar T)=
2\Lambda(\tau\wedge\bar T+\bar\tau\wedge T)-2i\tau\wedge\bar\tau
-i T\circ \bar T
-2iT\bar T.
\eea
We shall also need
\bea
\Lambda^3iT\wedge\bar T=6|\tau|^2-3|T|^2.
\eea

\subsection{Evaluation of the coefficients $A$ and $B$}

It is now easy to assemble all the terms and arrive at the final formula $B=0$ and $A$ is given by
\be
A=i\tilde Ric(\eta)+{i\over 2}T\circ\bar{T}.
\ee
We note that a simpler version of some of these identities when $m=3$ appeared in \cite{FP} and was instrumental in the proof of an upper bound for $\|\Omega\|_\o$.
Altogether the evolution equation for $\eta$ is
\be
\p_t\eta_{\bar kj}=-{1\over m-1}\left(\tilde{R}_{\bar kj}+{1\over 2}(T\circ\bar{T})_{\bar kj}\right).
\ee
This is the flow (\ref{eta-flow}) stated in Theorem \ref{main}. Q.E.D.

\section{Proof of Theorem \ref{unification}}
\setcounter{equation}{0}

Part (i) follows immediately from the fact that the Chern-Ricci tensor $\tilde R_{\bar kj}$ can be expressed in local coordinates as
\bea
\tilde R_{\bar kj}=-g^{p\bar q}\p_p\p_{\bar q} g_{\bar kj}+\cdots
\eea
where $\cdots$ denote terms with fewer derivatives. Part (ii) follows from the fact that, if $\tilde\eta(0)$ is K\"ahler, then the K\"ahler-Ricci flow $\p_t\tilde\eta=-R_{\bar kj}(\tilde\eta)$ admits a solution $\tilde\eta(t)$ which is K\"ahler for any $t$ in a small time interval near $t=0$. Since $T(\tilde\eta)=0$ and $\tilde R_{\bar kj}(\tilde\eta)=R_{\bar kj}(\eta)$, the same $\tilde\eta(t)$ satisfies the flow
(\ref{eta-flow}) if we take the same initial data $\eta(0)=\tilde\eta(0)$. By uniqueness of the solution of the flow (\ref{eta-flow}), it follows that $\eta(t)=\tilde\eta(t)$ for all time, and $\eta(t)$ is a solution of the K\"ahler-Ricci flow, as claimed. Part (iii) follows in the same way, using now Theorem \ref{main}. Indeed, if $\eta(0)$ is conformally balanced in the sense of Theorem \ref{unification}, then the corresponding $\omega(0)$ is conformally balanced in the sense of \cite{PPZ2}.
By Theorem 1 of \cite{PPZ4}, the Anomaly flow (\ref{anomaly} for $\omega$ admits a unique smooth solution $\o(t)$ for some small time interval near $0$. By Theorem \ref{main}, the corresponding $\eta(t)$ is a smooth, conformally balanced solution to the flow (\ref{eta-flow}). By uniqueness, this solution coincides with the solution known to exist by parabolicity. In particular, the conformally balanced condition is preserved for all $\eta(t)$. Q.E.D.

\section{Remarks}
\setcounter{equation}{0}

It may be interesting to find another, more direct, proof of Part (iii) of Theorem \ref{unification}, namely that the flow (\ref{eta-flow}) preserves the conformally balanced condition, instead of appealing to Theorem \ref{main} and the uniqueness of solutions. This does not appear evident, although it can for example be done for Part (ii). By deriving the flow for $|T|^2$ and applying the maximum principle, we can indeed show directly that the K\"ahler property is preserved. We reproduce the key calculations below, as the flow of the torsion is crucial in non-K\"ahler geometry, and the resulting formulas may be useful in other contexts. They are also comparatively simpler than the formulas for the flow of the torsion derived in \cite{PPZ2} under the conformally balanced condition.

\subsection{The flow of the torsion}

Consider then the flow
(\ref{eta-flow}), for general metrics $\eta$, not necessarily K\"ahler or conformally balanced. Introduce the notation $\eta=ig_{\bar kj}dz^j\wedge d\bar z^k$, and
write the flow (\ref{eta-flow}) as
\bea
\p_t\eta=-{1\over m-1}i(\tilde Ric(\eta)+{1\over 2}(T\circ\bar T))
\eea
Since $T=i\p\eta$, this implies immediately
\bea
\label{eta-flow-T}
\p_tT={1\over m-1}(\p\tilde Ric+{1\over 2}\p(T\circ\bar T))
\eea
For general Hermitian metrics, we have the following Bianchi identity
\bea
R_{\bar\ell m\bar kj}=R_{\bar\ell j\bar km}+\na_{\bar\ell}T_{\bar kjm}
=
R_{\bar kj\bar\ell m}+\na_j\bar T_{m\bar k\bar\ell}+\na_{\bar\ell}T_{\bar kjm}
\eea
and hence
\bea
\tilde R_{\bar kj}=R_{\bar kj}-\na_j\bar\tau_{\bar k}+\na^mT_{\bar kjm}
\eea
By Bochner-Kodaira formulas (see Appendix D), we have
\bea
\p^\dagger T_{\bar kj}=-\na^mT_{\bar kjm}+\bar\tau^mT_{\bar kjm}-{1\over 2}(T\circ\bar T)_{\bar kj}
\eea
Since the form $Ric(\eta)$ is closed and the form $\na_j\bar\tau_k$ is $\p$-exact, the right hand side of the equation (\ref{eta-flow-T}) can be expressed in terms of $\p^\dagger T$ and the differentials of $\bar\tau^mT_{\bar kjm}$ and $(T\circ\bar T)_{\bar kj}$ alone. We find that (\ref{eta-flow-T}) becomes
\bea
\label{eta-flow-T1}
\p_t T={1\over m-1}(-\p\p^\dagger T+\p(\bar\tau\cdot T)) \label{T}
\eea
where we have introduced the notation $\bar\tau\cdot T$ for the $(1,1)$-form defined by
\bea
(\bar\tau\cdot T)_{\bar\alpha\beta}=\bar\tau^\gamma T_{\bar\alpha\beta\gamma}.
\eea
Since $T$ is $\p$-closed, the operator $\p\p^\dagger$ on $T$ can be equated with the Hodge Laplacian
$\Box=\p\p^\dagger+\p^\dagger\p$ on $(2,1)$-forms, so we have obtained a parabolic diffusion equation for $T$.

\medskip
To show that the K\"ahler condition is preserved by the flow (\ref{eta-flow}), we need the evolution of $|T|^2$. For this, we again use a Bochner-Kodaira formula to convert the Hodge-Laplacian $\Box$ into the Laplacian $\Delta_c=g^{p\bar q}\na_{\bar q}\na_p$, modulo lower order terms,
\bea
(\p\p^\dagger T)_{\bar kjm}&=&\nabla_m(\p^\dagger T)_{\bar kj}-\nabla_j(\p^\dagger T)_{\bar km}+T^s_{~mj}(\p^\dagger T)_{\bar ks} \nonumber\\
&=&-\Delta_c T_{\bar kjm}+T_{\bar kjl}\tilde{R}^l_{~m}-T_{\bar kml}\tilde{R}^l_{~j}+g^{s\bar t}(R_{\bar tl\bar km}T^l_{~sj}-R_{\bar tl\bar kj}T^l_{~sm}) \nonumber\\
&&+\p(\bar\tau\cdot T)_{\bar kjm}+\frac{1}{2}\left(T^l_{~jm}(T\circ\bar T)_{\bar kl}+g^{s\bar t}g^{p\bar q}(\bar T_{m\bar t\bar q}\nabla_jT_{\bar ksp}-\bar T_{j\bar t\bar q}\nabla_mT_{\bar ksp})\right) \nonumber\\
&&+g^{s\bar t}g^{p\bar q}T_{\bar ksp}(R_{\bar qj\bar tm}-R_{\bar qm\bar tj}). \nonumber
\eea
Using the flow (\ref{eta-flow-T1}), we find
\bea
((m-1)\p_t-\Delta_c)|T|^2&=&-|\nabla T|^2-|\bar\nabla T|^2-\frac{1}{2}|T\circ\bar T|^2+\langle T\circ\bar T,T\bar T\rangle
\\
&&-2\Re(\bar{T}^{\bar kjm}(g^{s\bar t}g^{p\bar q}\bar{T}_{m\bar t\bar q}\nabla_jT_{\bar ksp}+2T^s_{~pj}R^{p}_{~s\bar km}+2T_{\bar ksp}R_{~j~m}^{p~s})).
\nonumber
\eea
In particular,
\bea
((m-1)\p_t-\Delta_c)|T|^2\leq C|T|^2(|T|^2+|Rm|).
\eea
for some constant $C$. The maximum principle implies that $T\equiv0$ if initially $T_0=0$, i.e., the K\"ahler condition is preserved.

We observe that it is easy to derive from the flow of $T$ the flows of $\tau$ as well as of the primitive component of $T$. For example, we have
\bea
(m-1)\p_t\tau_j&=&-\Box\tau_j+\nabla_j\left(|\tau|^2+\frac{1}{2}|T|^2\right)+T^s_{~pj}R'^p_{~s} \nonumber\\
&&+g^{p\bar k}\bar\tau_{\bar k}(\nabla_p\tau_j-\nabla_j\tau_p+T^s_{~pj}\tau_s).
\eea
In the conformally balanced case, we have $\p\tau=0$ and $R'^p_s=0$, and this flow reduces to
\bea
(m-1)\p_t\tau_j=-\Box\tau_j+\nabla_j\left(|\tau|^2+\frac{1}{2}|T|^2\right).
\eea
This results in the following flow for $|\tau|^2$,
\bea
((m-1)\p_t-\Delta_c)|\tau|^2&=&-|\nabla\tau|^2-|Ric|^2+\langle\tilde{R}ic-Ric+\frac{1}{2}T\circ\bar T,i\tau\wedge\bar\tau\rangle \nonumber\\
&&+\frac{1}{2}\langle \tau,\nabla|T|^2\rangle+\frac{1}{2}\overline{\langle \tau,\nabla|T|^2\rangle}+g^{j\bar k}(\bar\tau_{\bar k}T^s_{~lj}R^l_{~s}+\tau_j\bar T_{~\bar p\bar k}^{\bar t}R_{\bar t}^{~\bar p})\nonumber
\eea
The right hand side can only be bounded above by $O(|\tau|)$, which indicates that, unlike the vanishing of the full torsion $T$, the vanishing of $\tau$ is not preserved. Indeed one can verify that the balanced condition $\tau=0$ is not preserved by running the Anomaly flow on generalized Calabi-Gray manifolds \cite{FHP} with balanced initial data.

\subsection{Shi-type estimates and long-time existence of the flow}

Shi-type estimates for the original Anomaly flow were derived in \cite{PPZ2}. For the present version (\ref{eta-flow}) in terms of the rescaled metric $\eta$, they become simpler to establish, and the same arguments as in \cite{PPZ2}, or the general results in \cite{ST}, imply the following statement: the flow in $\eta$ will continue to exist, unless there is a time $T>0$ and a sequence $(z_j,t_j)$ with $t_j\to T$, with
\bea
|Rm(z_j)|^2_{\eta(t_j)}+|T(z_j)|^2_{\eta(t_j)}+|\na T(z_j)|^2_{\eta(t_j)}\to\infty.
\eea
Now under a conformal change of metrics $\o_f=e^f\o$, the torsion and curvature transform as follows
\bea
&&
T^\ell{}_{jk}(\o_f)=T^\ell{}_{jk}(\o)+f_j\delta^\ell{}_k-f_k\delta^\ell{}_j, \nonumber\\
&&
R_{\bar kj\bar pq}(\o_f)=e^f(R_{\bar kj\bar pq}(\o)+f_{\bar kj}\o_{\bar pq}).
\eea
In the case at hand, $\eta=\|\Omega\|_\o\o$, so it is easy to work out the previous conditions, and find that the flow will continue to exist unless there is a time $T$ and a sequence $(z_j,t_j)$ with $t_j\to T$ satisfying
\bea
{|Rm(z_j)|^2_{\omega(t_j)}
\over\|\Omega(z_j)\|^2_{\omega(t_j)}}+{|T(z_j)|^2_{\omega(t_j)}\over\|\Omega(z_j)\|_{\omega(t_j)}}+{|\na T(z_j)|^2_{\omega(t_j)}\over\|\Omega(z_j)\|^2_{\omega(t_j)}}\to\infty.
\eea
This is a more succinct, and perhaps more natural formulation of the criterion for the appearance of singularities found in \cite{PPZ2}, which involved the four quantities $|Rm(z_j)|_{\omega(t_j)}$, $|T(z_j)|_{\omega(t_j)}$, $|\na T(z_j)|_{\omega(t_j)}$, and $\|\Omega(z_j)\|_{\omega(t_j)}$ separately.

\subsection{Two questions by Ustinovskiy}

In Ustinovskiy's thesis \cite{U}, he raised two questions (Question 6.15 and Problem 6.16) about periodic solutions and stationary points of Hermitian curvature flow. As a result of our Theorems 1 and 2, we can answer these questions.

\begin{proposition}
All periodic solutions are stationary, which are Ricci-flat K\"ahler metrics.
\end{proposition}

The proof goes as follows: Suppose that we have a periodic solution. It follows from the results of Ustinovskiy and Theorems 1 and 2 that the flow is exactly the Anomaly flow with $\alpha'=0$ and with conformally balanced initial data. Therefore we have monotonicity formulae as introduced in \cite{FP}. Periodic solutions imply that all the monotone quantities are actually constants, which in turn gives us an equation from the monotonicity formula. This equation can only be satisfied by Ricci-flat K\"ahler metrics, which are the only stationary points of the Anomaly flow. Note that the second part of this proposition has been established in several different ways in the literature, including by an integration by parts, and by many authors including \cite{CHSW},
\cite{FT}, \cite{MT}, and \cite{PPZ4}.

\begin{appendix}

\section{Conventions and preliminaries}
\setcounter{equation}{0}

If $\eta$ is a Hermitian metric on $X$, its curvature is defined by $R_{\bar kj}{}^p{}_q=-\p_{\bar k}(g^{p\bar m}\p_jg_{\bar m q})$.
Its Ricci curvature $R_{\bar kj}(\eta)$ is defined by
\bea
R_{\bar kj}(\eta)=-\p_j\p_{\bar k}\log\,\eta^m=\p_j\p_{\bar k}\log\,\|\Omega\|_\eta^2
\eea
and the Ricci form $Ric(\eta)$ is defined by
\bea
Ric(\eta)=R_{\bar kj}(\eta)dz^j\wedge d\bar z^k
=\p_j\p_{\bar k}\log\,\|\Omega\|_\eta^2 dz^j\wedge d\bar z^k=
\p\bar\p \log\,\|\Omega\|_\eta^2.
\eea
As in \cite{PPZ2}, the other notions of Ricci curvature are defined by $\tilde R_{\bar kj}=R^p{}_{p\bar kj}$, $R_{\bar kj}'=R_{\bar k}{}^p{}_{pj}$, $R_{\bar kj}''=R^p{}_{j\bar kp}$.
Our conventions for the torsion of $\eta$ are
\bea
T=i\p\eta\equiv{1\over 2}T_{\bar kjm}dz^m\wedge dz^j\wedge d\bar z^m
\eea
In particular $T^m{}_{jp}=\Gamma^m_{jp}-\Gamma^m_{pj}$, with $\Gamma^m_{jp}=g^{m\bar q}\p_jg_{\bar qp}$.
We set
\bea
\tau_l=g^{j\bar k}T_{\bar kjl}.
\eea
The norm $\|\Omega\|_\o^2$ with respect to a given Hermitian metric is defined as usual as
\bea
i^{m^2}\Omega\wedge\bar\Omega=\|\Omega\|_\o^2{\o^m\over m!}.
\eea

\section{The Hodge operator $\Lambda$}
\setcounter{equation}{0}

We define the operator $\Lambda^q$ from $(p,p)$-forms to $(p-q,p-q)$-forms by
\bea
(\Lambda^q\Phi)_{\bar j_1k_1\cdots\bar j_{p-q} k_{p-q}}
=i^{-q}\prod_{\alpha=p-q+1}^p g^{k_\alpha\bar j_\alpha}\Phi_{\bar j_1k_1\cdots\bar j_pk_p}
\eea
for
\bea
\Phi={1\over (p!)^2}\Phi_{\bar j_1k_1\cdots\bar j_pk_p}dz^{k_p}\wedge dz^{\bar j_p}\wedge\cdots dz^{k_1}\wedge d\bar z^{j_1}
\eea
Note that $\Lambda$ maps real forms to real forms, and that $\Lambda^p\Phi=\<\Phi,\eta^p\>={\rm Tr}\,\Phi$. In terms of $\Lambda$, the previous torsion component $\tau$ can be written as $\tau=i\Lambda T$.

\section{The Hodge $\star$ operator}
\setcounter{equation}{0}

Let $\alpha$, $\Phi$, and $\Psi$ be $(1,1)$-forms, $(2,2)$-forms, and $(3,3)$-forms respectively. Then we have the following identities (the detailed derivations can be found in \cite{PPZ4})
\bea
\star(\alpha\wedge\eta^{m-2})&=& (m-2)! (-\alpha+(\Lambda\alpha)\eta)
\nonumber\\
\star(\Phi\wedge\eta^{m-3})&=&(m-3)!(-\Lambda\Phi+{1\over 2}(\Lambda^2\Phi)\eta)
\nonumber\\
\star(\Psi\wedge\eta^{m-4})&=&(m-4)!(-{1\over 2}\Lambda^2\Psi+{1\over 6}(\Lambda^3\Psi)\eta)
\eea

Let $\tau$ and $T$ be $(1,0)$-forms and (2,1)-forms respectively. Then we also have
\bea
\star(\tau\wedge\eta^{m-2})&=&-i(m-2)!\tau\wedge\eta
\nonumber\\
\star(T\wedge\eta^{m-3})&=&i(m-3)!(-\Lambda T\wedge\eta+T)
\eea

\section{The operator $\p^\dagger$ and Bochner-Kodaira formulas}
\setcounter{equation}{0}

First we work out the operator $\p^\dagger$ on various spaces of forms. The basic formula is the following integration by parts formula for general Hermitian metrics
\bea
\int_X \nabla_jV^j \o^n
=
\int_X\tau_j V^j\o^n
\eea
where $V^j$ is a vector field.

\medskip

To get e.g. the $\p^\dagger$ on $(1,0)$-forms, we take $V^j=f\overline{\psi_k}g^{j\bar k}$ where $f$ and $\psi_k$ are respectively an arbitrary scalar function and an arbitrary $(1,0)$-form. Them
\bea
\int_X\nabla_j(f\overline{\psi_k}g^{j\bar k})\o^n=\int_X \tau_jf\tau_j\overline{\psi_k}g^{j\bar k}o^n
\eea
which can be rewritten as
\bea
\int_X (\nabla_j f)\overline{\psi_{\bar k}}g^{j\bar k}\o^n
+
\int_X f\overline{g^{k\bar j}\nabla_{\bar j}\psi_{k}}\o^n=\int_X f\overline{g^{k\bar j}\bar\tau_{\bar j}\psi_k}\o^n
\eea
This means that
\bea
\p^\dagger\psi=-g^{k\bar j}\na_{\bar j}\psi_{k}+g^{k\bar j}\bar\tau_{\bar j}\psi_k
\eea
for any $(1,0)$-form $\psi$.
More generally, we have

\begin{lemma}
Suppose $\alpha$ is a $(p,q)$-form, then
\bea
(\p\alpha)_{\bar t_q\dots\bar t_1s_{p+1}\dots s_1}&=&\sum_{k=1}^{p+1}(-1)^{k-1}\p_{s_k}\alpha_{\bar t_q\dots\bar t_1s_{p+1}\dots\widehat{s_k}\dots s_1} \nonumber\\
&=&\sum_{k=1}^{p+1}(-1)^{k-1}\nabla_{s_k}\alpha_{\bar t_q\dots\bar t_1s_{p+1}\dots\widehat{s_k}\dots s_1}+\sum_{l<k}(-1)^{k}T^s_{~s_ls_k}\alpha_{\bar t_q\dots\bar t_1s_{p+1}\dots\widehat{s_k}\dots s_{l+1}ss_{l-1}\dots s_1}.\nonumber
\eea
and
\bea
(\p^\dagger\alpha)_{\bar t_q\dots\bar t_1s_{p-1}\dots s_1}&=&-g^{s\bar k}\nabla_{\bar k}\alpha_{\bar t_q\dots\bar t_1s_{p-1}\dots s_1s}+g^{s\bar k}\bar\tau_{\bar k}\alpha_{\bar t_q\dots\bar t_1s_{p-1}\dots s_1s} \nonumber\\
&&+\frac{1}{2}\sum_{l=1}^{p-1}(-1)^l\bar{T}_{s_l}^{~~cd}\alpha_{\bar t_q\dots\bar t_1cds_{p-1}\dots \widehat{s_{l}}\dots s_1}.\nonumber
\eea
\end{lemma}

For example, we have:

$\bullet$ If $\alpha$ is a (1,0)-form, then
\bea
(\p\alpha)_{pq}&=&\nabla_q\alpha_p-\nabla_p\alpha_q+T^s_{~qp}\alpha_s, \nonumber\\
\p^\dagger\alpha&=&-g^{p\bar k}\nabla_{\bar k}\alpha_p+g^{p\bar k}\bar{\tau}_{\bar k}\alpha_p.\nonumber
\eea

$\bullet$ If $\beta$ is a (2,0)-form, then
\bea
(\p^\dagger\beta)_l&=&-g^{p\bar k}\nabla_{\bar k}\beta_{lp}+g^{p\bar k}\bar{\tau}_{\bar k}\beta_{lp}-{1\over 2}\bar{T}_{l\bar a\bar b}\beta_{cd}g^{c\bar a}g^{d\bar b}. \nonumber
\eea

$\bullet$ If $\psi$ is a $(2,1)$-form, then
\bea
(\p^\dagger\psi)_{\bar\alpha\beta}
=
-g^{\gamma\bar j}\na_{\bar j}\psi_{\bar\alpha\beta\gamma}
+
g^{\gamma\bar j}\bar\tau_{\bar j}\psi_{\bar\alpha\beta\gamma}
-
{1\over 2}
\bar T_{\beta\bar j\bar m}\psi_{\bar\alpha\gamma\delta}g^{\gamma\bar j}g^{\delta\bar m}
\eea

\end{appendix}

\bigskip
\noindent
{\bf Acknowledgements} The authors would like to thank Sebastien Picard and Xiangwen Zhang for sharing generously their insights, and for joint work on many projects closely related to this paper. The second-named author would like to thank Carlo Maccaferri, Rodolfo Russo, and the Galileo Galilei Institute for Theoretical Physics in Florence for their kind hospitality where part of this work was done.

\end{document}